\documentclass[11pt,twoside]{article}
\usepackage{amsmath}
\usepackage{amssymb}
\usepackage{graphicx}
\usepackage{enumerate}
\usepackage[latin1]{inputenc}
\usepackage{latexsym}
\usepackage{graphicx,epsfig}
\input epsf

\def\limiten{\renewcommand{\arraystretch}{0.5}
\begin{array}[t]{c}\stackrel{}{\longrightarrow} \\
{\scriptstyle n\rightarrow
\infty}\end{array}\renewcommand{\arraystretch}{1}}

\def\limiteloin{\renewcommand{\arraystretch}{0.5}
\begin{array}[t]{c}\stackrel{{\cal D}}{\longrightarrow} \\
{\scriptstyle n\rightarrow
\infty}\end{array}\renewcommand{\arraystretch}{1}}

\newcommand{\be}{\begin{equation}}
\newcommand{\ee}{\end{equation}}
\newcommand{\bd}{\begin{displaymath}}
\newcommand{\ed}{\end{displaymath}}

\newcommand{\ba}{\begin{eqnarray}}
\newcommand{\ea}{\end{eqnarray}}
\newcommand{\ban}{\begin{eqnarray*}}
\newcommand{\ean}{\end{eqnarray*}}

\newcommand{\ds}{\displaystyle}
\newcommand{\LL} {I\!\!L}
\newcommand{\R} {I\!\!R}
\newcommand{\E} {I\!\! E}
\newcommand{\N} {I\!\! N}

\renewcommand{\arraystretch}{.8}

\renewcommand{\Box}{\hfill\rule{0.25cm}{0.25cm}} 

\newtheorem{Prop}{Proposition}[section]
\newtheorem{lem}{Lemme}[section]
\newtheorem{Theo}{Théorème}[section]

\newtheorem{rem}{Remarque}[section]

\newenvironment{dem}{\ \\ {\bf Démonstration. }}
{\Box\par\medskip\noindent}

\def\1{{\bf 1}}

\begin{document}
\title{\bf Estimation non-paramétrique de la densité spectrale d'un processus gaussien échantillonné aléatoirement}
\date{}
\maketitle \vspace{-0.5cm}
\begin{center}
Jean-Marc~BARDET${}^*$,\; Pierre, R.~BERTRAND${}^{**}$ and Véronique~BILLAT${}^{\dag}$
\end{center}
${}^*$  {\it SAMOS-MATISSE - UMR CNRS 8595, Universit\'e
Panth\'eon-Sorbonne (Paris I),
90 rue de Tolbiac, 75013 Paris Cedex, France, E-mail: bardet@univ-paris1.fr} \\ \\
${}^{**}$ {\it Laboratoire de Math\'ematiques - UMR CNRS 6620,
Universit\'e Blaise Pascal (Clermont-Ferrand II), 24 Avenue des
Landais, 63117 Aubi\`ere Cedex, France. E-mail:
Pierre.Bertrand@math.univ-bpclermont.fr}\\ \\
${}^{\dag}$ {\it UBIAE Laboratoire  INSERM 902 et
Université d'\'Evry.} 
 \pagestyle{myheadings}
\markboth{Estimation non paramétrique de la densité spectrale}{
J-M. Bardet,\; P.R. Bertrand and V. Billat} ~\\ ~\\
{\bf Résumé:}\\
En utilisant une analyse en ondelette, on construit un
estimateur non-paramétrique de la densité spectrale d'un processus
gaussien à accroissements stationnaires. Dans un premier temps, on
considère le cas ``idéal'' de l'observation d'une trajectoire
en temps continu et on donne un théorème de la limite central ponctuel et  une estimation de l'erreur quadratique moyenne  intégrée.  Puis, afin de mieux correspondre aux applications, on construit un
second estimateur à partir de l'observation d'une trajectoire du
processus à des instants discrets, aléatoires et irrégulièrement
espacés. Cet estimateur est obtenu en remplaçant les coefficients
d'ondelette par leurs discrétisations. Nous donnons un second théorème de la limite central,
avec une vitesse différente. Des simulations et une application
à la modélisation du rythme cardiaque des marathoniens sont
 présentées.
\\
~\\
 {\bf Abstract (english):}
\\
From a wavelet analysis, one derives a nonparametrical estimator
for the spectral density of a Gaussian process with stationary
increments. First, the idealistic case of a continuous time path of the process is considered. A
punctual Central Limit Theorem (CLT) and an estimation of the Mean Integrate Square Error (MISE) are established. Next, to
fit the applications, one considers the case where one observes a path at random times.
One built a second estimator obtained by replacing the wavelet coefficients by their
discretizations. A second CLT and the corresponding estimation of the
MISE are provided. Finally, simulation results  and an application on the heartbeat time series of marathon runners
are presented.

\section{Introduction}
Soit ${\bf X}= \{X(t), t\in\R\}$ un processus à temps continu
gaussien centré à accroissements stationnaires
Il admet une représentation (voir Yaglom, 1958)
\be
\label{repr:harmonizable} X(t)= \int_{\R}
\big(e^{it\xi} - 1\big)\cdot f^{1/2}(\xi) \,
dW(\xi),~~~~\mbox{pour tout}~~t \in \R, \ee
où $W(dx)$ est une mesure de Wiener
complexe adaptée afin que $X(t)$ soit réel et $f$ est une fonction paire, positive appelée densité spectrale.
On observe une trajectoire $(X_{t_0},\ldots,X_{t_n})$ du processus $X$ à des instants $t_i$, $i=0,\dots,
n$ irrégulièrement espacés ou aléatoires. Le but de cet article est de proposer un estimateur de
la densité spectrale $f(\xi)$ pour une fréquence $0<\xi<\infty$ et d'étudier sa convergence.

A titre d'exemple paradigmatique, si $f(\xi)= \sigma^2\,|\xi|^{-(2H+1)}$ pour tout $\xi \in \R^*$, alors le processus $ X$
est un mouvement brownien fractionnaire (mbf) d'indice de Hurst
$H$ et de paramètre d'échelle $\sigma$; la fonction de densité
spectrale est entièrement déterminée par les deux seuls
paramètres $(H, \,\sigma)$. On se trouve dans un cadre
d'estimation paramétrique de la densité spectrale $f$. De plus, la
fonction $\ln f(\xi)$ est une fonction affine de $\ln(|\xi|)$, ce
qui correspond à l'autosimilarité du processus.

L'autosimilarité est très souvent supposée sans vraiment être
statistiquement vérifiée ou théoriquement justifiée.
Dans cet article, nous nous intéressons à l'estimation de la
fonction $f(\xi)$ dans un cadre non-paramétrique. Cette approche peut
permettre ensuite de vérifier graphiquement l'autosimilarité, ou, par exemple,
de déterminer une ou plusieurs bandes de fréquences de fréquences
$]\omega_i,\omega_{i+1}[$ telle que $f(\xi) =
\sigma_i^2\cdot|\xi|^{-(2H_i+1)}$ pour $\omega_i<\xi<\omega_{i+1}$.

Le plan de cet article est le suivant. Dans la section 2, nous précisons le
modèle statistique et les hypothèses utilisées. Dans la
section 3, on montre un théorème de la limite central (TLC) pour
l'estimation de la fonction de densité spectrale à l'aide d'une
analyse en ondelette théorique, c'est-à-dire avec des coefficients
d'ondelette calculés par une intégrale sur $\R$. Dans la section
suivante, nous montrons que la discrétisation des coefficients
d'ondelette découlant de l'observation à des instants discrets
aléatoires conduit à établir un second TLC. Dans la dernière section, nous donnons des
exemples numériques.
\section{Hypothèses et notations}
Soit $X$ un processus gaussien centré à accroissements
stationnaires défini par l'équation (\ref{repr:harmonizable}). On
suppose que sa densité spectrale $f$ vérifie les conditions (F1)
et (F2) ci-dessous:
\begin{description}
\item[(F1) ] {\em La fonction $f$ est borélienne, paire positive et
vérifie}
\be \label{cond:f} \int_{\R} \big(1\wedge
|\xi|^2\big)\cdot f (\xi) \, d \xi < \infty. \ee
\item[(F2) ]
{\em Sur $\R^+$ la fonction $f$ est càdlàg et $\mathcal{C}^1$ en
dehors d'un nombre fini de points  $0<\omega_1<\dots<\omega_K$
avec $K
\in
\N$, sa derivée $f'$ vérifie la condition
 \be
\label{cond:f:prime} \int_{\R} \left(1\wedge |\xi|^3\right)\cdot
\left|f'(\xi)\right| \, d \xi < \infty. \ee
De plus, on a un contrôle du comportement à haute fréquence: il
 existe deux constantes $C_0,C_0'>0$, une fréquence $\omega_K>0$ et un réel $H
>0$, tel que pour $|\xi|\ge \omega_K$ }
\ba \label{majoration:f:infini} f(\xi)\,\le\, C_0\,
|\xi|^{-(2H+1)}&\quad\mathrm{et}\quad& f'(x)\,\le\, C_0'\
|\xi|^{-(2H+2)}.\ea
\end{description}
On dispose de l'observation du processus à des instants discrets
sur un intervalle $[0, T_{n}]$,  {\it
i.e.} \\ $\ds \big (X(t_0),X(t_1),\ldots,X(t_n)\big )$ est connue
avec $t_0 =0<t_1 <\dots<t_{n} =T_n$. Les instants d'observations
$t_i$ sont aléatoires et irrégulièrement espacés. On aura besoin
que l'hypothèse (S$(s)$) ci-dessous soit vérifiée pour un réel $s\in
[1,\infty]$.
\begin{description}
\item[(S$(s)$)] {\em
Il existe une suite de nombres réels strictement positifs
$(\delta_n)_{n \in \N}$ telle que $\ds\delta_n \limiten 0$ et une suite de
variables aléatoires positives identiquement distribuées et indépendantes $(L_k)_{k\in
\N }$ telle que
$$\ds t_{k+1}-t_{k} :=
\delta_n \, L_k\quad\mbox{et $\E \, L_k=1$ pour $k \in
\{0,1,\ldots,n-1\}$}.
$$
De plus, les v.a. $(L_k)_{k \in \N}$ sont indépendantes de la
tribu ${\cal F}_X$  engendrée par le processus avec  ${\cal F}_X
:=
\sigma
\big
\{X(t),\, t
\in
\R \big \}$ et il existe deux constantes $0<m_s<M_s<\infty$ telles que
$$
m_s\leq \|L_{k}\|_s \leq M_s~~~\mbox{pour tout  $k \in \N$}.
$$
où  pour une v.a. $Z$ et un réel $\alpha\in (0,\infty]$, on
note$\| Z\|_\alpha:=\big ( \E (|Z|^\alpha)
\big ) ^{1/\alpha}$ quand   $\E (|Z|^\alpha) <\infty$ et
$\| Z\|_\infty=\sup_{\omega
\in
\Omega}|Z(\omega)|$ par convention.}
\end{description}
Un estimateur de la densité spectrale $f(\xi)$ pour une fréquence
$0<\xi <\infty$ est construit à partir de l'analyse en ondelette
du processus $X$. Soit $\psi$ une {\em "ondelette mère"}, et pour n'importe quels échelle et
décalage $(a,b)\in
\R_+^*\times
\R$, on définit respectivement les coefficients d'ondelette
``théoriques'' (pour une trajectoire à temps continu) et ``empiriques'' (pour une trajectoire à temps discrets) par
\ba
\label{def:coeff:W}
d_\psi(a,b&:=&  a^{-1/2} \int_{\R} \psi\left(
\frac{t-b}{a}\right)\,X(t)\,dt\\
\label{def:coeff:W:discrets}
\mathrm{et}\hspace{3cm} e_\psi(a,b) &:=&  a^{-1/2} \sum
_{i=0}^{n-1}\left(\int_{t_i}^{t_{i+1}}  \psi
\Big(\frac {t -b} a \Big)\, dt\right) \, X \big(t_i\big).\qquad\qquad\qquad
\ea
Suivant les différents cas examinés par la suite, l'ondelette $\psi$ vérifiera des hypothèses issues de la famille suivante~:
\begin{description}
\item[W$(m,q,r)$] {\em $\psi:~\R
\mapsto \R$ est une fonction vérifiant les trois conditions:
\begin{itemize}
\item pour tout entier $n\le m$, $\displaystyle{\int_{\R} \left |t^n\psi(t)\right |dt
<\infty}$, $\displaystyle{\int_{\R}  t^n\psi(t)\,dt =0}$ et
$\displaystyle{\int_{\R} \left |t^{m+1}\,\psi(t)\right |dt
<\infty}$;
\item il existe une constante $C_{\psi} >0$ telle que pour tout $t\in \R$
\ban \big(1+|t|\big)^q \cdot
 \big| \psi(t)\big|
&\le& C_{\psi}.
 \ean
\item il existe une constante $C'_{\psi} >0$ telle que pour tout $\xi\in \R$
\ban \big(1+|\xi|\big)^r \cdot
\big(\big|\widehat{\psi}(\xi)\big|+\big|\widehat{\psi}'(\xi)\big|\big)
&\le& C'_{\psi}.
 \ean
\end{itemize}}
\end{description}
Cette famille d'hypothèses est également  utilisée dans Bardet et
Bertrand (2008). Dans cette présentation, nous ferons deux
hypothèses supplémentaires qui permettent de simplifier les énoncés et les
démonstrations~:
\begin{description}
\item[(T)] {\em
\begin{itemize}
\item La durée d'observation $T_n$ est déterministe.
\item
La transformée de Fourier de l'ondelette mère, $\widehat{\psi}$ a
un support compact.
\end{itemize}}
\end{description}
\begin{rem}
Les différentes hypothèses imposées sont faibles. 
\begin{itemize}
  \item L'hypothèse (F1) est
équivalente à l'existence de l'intégrale stochastique
(\ref{repr:harmonizable}), l'hypothèse (F2) autorise à la fois des
mbf, des processus à bande finie (voir Rachdi, 2004), les processus
localement fractionnaires utilisés par Kammoun {\it et al.} (2007) pour
modéliser les battements cardiaques ou des processus fractionnaires
multi-échelles définis dans Bertrand et Bardet (2001). 
  \item L'hypothèse S$(s)$ est satisfaite pour tout
$s<\infty$ dans le cas où les intervalles d'inter-arrivée suivent une loi
exponentielle (on parlera alors d'échantillonnage poissonnien) et pour $s=\infty$ dans le cas déterministe ou le cas
intervalles où la loi des intervalles d'inter-arrivée a un support
inclus dans $[a\delta_n,b\delta_n]$ avec $0<a<b$; la condition $\ds\delta_n
\limiten 0$ est naturelle. Seule l'hypothèse d'indépendance
des instants d'observations par rapport au processus est
restrictive pour les applications en finance, mais elle est
toujours faite, nous renvoyons sur ce point à Aït-Sahalia
et Mykland (2008). A notre connaissance, l'hypothèse S(s) avec $s<\infty$ est nouvelle
et tous les travaux sur les mbf ou les  diffusions observés à des
instants irréguliers ou  aléatoires  supposent vérifiée
S($\infty$), voir par exemple Begyn (2005), Blanke et Vial (2007)
ou Hayashi et Yoshida (2005).
  \item  L'hypothèse W$(m,q,r)$ est
classique et peu contraignante, voir la discussion ci-dessous.
\item
L'hypothèse (T) est purement technique et peut être supprimée au
prix de notations plus lourdes. Supposer déterministe  la durée
d'observation $T_n$ est standard. Cela n'empêche
pas les temps d'inter-arrivée d'être aléatoires comme par exemple dans le cas d'un processus ponctuel.
Dans les applications, que
ce soit en finance lors de l'observation des cours d'une action ou
en biologie lors de l'observation du rythme cardiaque d'un
individu sain ou malade, on a une durée d'observation $T_n$
déterministe ($8$h pour les cours de bourse, $1$h ou $24$h pour les
battements de coeur). On peut donner dependant un exemple où $T_n$ est aléatoire~: Kammoun {et al.} (2007)
ont étudié le
rythme cardiaque de coureurs de marathon, $T_n$ correspondant alors à la
durée du marathon, qui est indubitablement aléatoire\footnote{Un
match de tennis a lui aussi une durée aléatoire, par contre un
match de football a une durée à peu près déterministe.}.
\end{itemize}
\end{rem}
\section{Analyse en ondelette d'une trajectoire continue}
{\bf Qu'appelle-t-on exactement une {\em "ondelette mère"}?} Il y
a en fait un léger abus, car nous n'avons pas besoin que la
fonction $\psi$ et ses {\em "filles"} engendrent $\LL^2(\R)$. Nous utilisons
uniquement des hypothèses W$(m, q, r)$ qui sont faibles et signifient que $\psi$ est une
fonction bien localisée en temps et en fréquence avec un certain
nombre de moments nuls, comme l'est par exemple une ondelette de Lemarié-Meyer.
La première étape consiste à obtenir une représentation
fréquentielle des coefficients:
\begin{Prop} \label{pte:d(a)}
Soient $X$ un processus gaussien défini par
(\ref{repr:harmonizable}) avec une densité spectrale $f$ vérifiant
(F1) et $\psi$ une ondelette vérifiant  la condition W$(0,0,0)$,
alors pour tout $(a,b)\in
\R_+^*\times
\R$, le coefficient d'ondelette $d_{X}(a,b)$ est bien défini par (\ref{def:coeff:W}) et
admet la représentation fréquentielle
\begin{eqnarray} \label{repr:dX}
d_{X}(a, b) =\sqrt{a}\, \int_{\R} e^{ib\xi}\,
\overline{\widehat{\psi}}(a \xi)\, f^{1/2}(\xi)\,dW(\xi).
\end{eqnarray}
\end{Prop}
\begin{dem} La démonstration résulte de l'application du Théorème
de Fubini stochastique, licite grâce aux hypothèses. Elle est
donnée en détail, tout comme les autres démonstrations des résultats de
cet article, dans Bardet et Bertrand (2008).
\end{dem}
On en déduit que, pour une échelle $a$ fixée, le processus $\ds\big(d_\psi(a,b)\big)_{b\in
\R}$ est stationnaire, gaussien, centré, et a pour variance
$$\ds
\mathcal{I}_{\psi}(a) :=  \int_{\R} |\widehat{\psi}(x)|^2\, f(x/a)\,
du.$$
On  montre également que, pour $a$ fixé, le processus $\ds\big(d_\psi(a,b)\big)_{b\in
\R}$ n'est pas à longue mémoire (alors que le processus $X$ peut l'être). Tout estimateur de la variance de
$d_\psi(a,\cdot)$ fournit alors une estimation de $\ds
\mathcal{I}_{\psi}(a)$.
\subsubsection*{Estimation de la variance des coefficients
d'ondelette}
Soit $(N_k)_{k \in \N}$ une suite d'entiers telle que $\ds N_n \limiten \infty$. Soit
$n \in \N^*$. Pour une famille de décalages
 $\ds b_1<b_2<\dots <b_{N_n}$, on estime $\ds
\mathcal{I}_{\psi}(a)$ par la variance empirique
\ba\label{def:I:NPsi}   I_{N_n,\psi}(a) &:=& N_n^{-1}
\,
\sum_{k=1}^{N_n}  d^2_{\psi}(a,b_k).
\ea
\subsubsection*{Comment choisir la famille de décalages
$(b_k)_{k=1,\dots,N_n}$~?}
Dans le cas d'une observation à des
instants déterministes régulièrement espacés, l'usage consiste à
prendre comme famille de décalages les instants d'observation,
sans que jamais la question du choix ne soit posée. Nous n'avons
pas réussi pour l'instant à résoudre la question du choix optimal
de la famille de décalages. Pour simplifier les calculs et la
présentation, nous supposons dans la suite que les décalages $b_k$
sont régulièrement espacés. On montre dans la section suivante que
l'erreur de discrétisation est contrôlée si
$T_n^\rho<b_1<\dots<b_{N_n}<T_n-T_n^\rho$ avec $\rho \in (3/4,1)$.
Ceci conduit au choix suivant:
\begin{description}
\item[(B1)]{\em
On se fixe un réel $ \rho \in (3/4,1)$, on note $\tau_n := T_n-2
T_n^\rho$ et on définit la famille
de décalages par}
\ba
\label{def:bk}
b_k &:=& T_n^\rho + \frac{k-1}{N_n-1}\times \tau_n\qquad pour\; k=
1,\dots,N_n.
\ea
\end{description}
Remarquons que $\tau_n = \big(b_{N_n}-b_1\big)$. La proposition suivante
donne une condition
suffisante pour que l'estimateur $I_{N_n,\psi}(a)$ satisfasse un TLC.
\begin{Prop}\label{TCL:ponctuel}
Soient $X$ un processus gaussien défini par
(\ref{repr:harmonizable}) avec une densité spectrale $f$ vérifiant
(F1) et (F2), $\psi$ une ondelette vérifiant  la condition W$(1, 0,
1/2)$ et $(b_i)_{1\leq i\leq N_n}$ vérifiant(\ref{def:bk}). Alors si
$\ds \tau_n\limiten \infty,\quad N_n\limiten \infty\quad
\mathrm{et}\quad \tau_n/N_n
\limiten 0$,
 \ba
\label{TCL1} \sqrt{\tau_n}\,  \big (I_{N_n,\psi}(a)
- \mathcal{I}_{\psi}(a) \big )\limiteloin
\mathcal{N}\Big (0\, ,\, 4\pi \, a^2 \int_{\R}
 \big |
\widehat{\psi}(a x)\big |^4
\, f(x)^2\,dx\Big ).
\ea
\end{Prop}
\subsubsection*{Estimation ponctuelle de la densité spectrale}
Dans l'expression de $\ds
\mathcal{I}_{\psi}(a)$, la fonction $x \mapsto
|\widehat{\psi}(x)|^2$ joue le rôle de la fenêtre spectrale et
on estime $\mathcal{I}_{\psi}(a)$ qui est une moyenne de $f$ au
voisinage de $x=1/a$. Pour estimer $f$ ponctuellement, il suffit
de disposer d'une suite d'ondelettes $\psi_\lambda$ vérifiant W$(1,
1/2, 0)$ telles que leurs transformées de Fourier
$\widehat{\psi_\lambda}$ convergent vers la masse de Dirac en $1$.
On montre ainsi le lemme suivant:
\begin{lem}\label{lemDL}
Soit $\psi$ une ondelette vérifiant la condition W$(1, 0, 1/2)$ telle que le support de
sa transformée de Fourier soit inclus dans un compact $[-\Lambda, \Lambda]$ avec $\Lambda>0$. Soit
la famille
$(\psi_\lambda)_{\lambda \in \R_+^*}$ telle que pour tous $\lambda \in \R_+^*$ et $t\in \R$, $\ds
\psi_\lambda(t) \,=\,  \big(1/ {\sqrt{\lambda}}\big)\, e^{i\,t}\cdot
\psi( t/ \lambda)$. De plus,  si
la fonction  $f$ est $\mathcal{C}^2$ dans un voisinage du point
$1/a$, on a le développement limité  quand $\lambda\to
\infty$
\ba
\label{DL:I:psi}
\mathcal{I}_{\psi_\lambda}(a) &=& \|\psi\|^2_{L^2} \times \left[f(1/a) + f'(1/a)\cdot \frac 1 {a\lambda}
+ \mathcal{O}\left((\lambda)^{-2}\right)\right]
\ea
\end{lem}
\begin{dem}
On a $\ds\widehat{\psi}_\lambda(\xi) = \sqrt{\lambda}\cdot
\widehat{\psi}\big(\lambda(\xi-1)\big)$ pour $\xi \in \R$.
La première condition de W(1,0,1/2) signifie que
$\widehat{\psi_\lambda}(0) =
\widehat{\psi_\lambda}'(0) =0$ et $\ds \int_{\R} |t|^2 |\psi_\lambda(t)|\, dt
<\infty$, or $\widehat{\psi_\lambda}(0) =
\sqrt{\lambda}\cdot\widehat{\psi_\lambda}(-\lambda) = 0$ et $\widehat{\psi_\lambda}'(0) =
 \lambda^{3/2}\cdot\widehat{\psi_\lambda}'(-\lambda)= 0$ pour $\lambda> \Lambda$, les autres conditions se vérifiant aisément.
D'autre part
$$
\mathcal{I}_{\psi_\lambda}(a) = \int_{\R}\lambda \,|\widehat{\psi}\big(\lambda(x-1)\big)|^2\, f(x/a)\,
dx =
\int_{\R} |\widehat{\psi}\big(u\big)|^2\times f\left(\frac 1 {a} + \frac u {\lambda a}\right)\,
du,
$$
puis, en utilisant la formule de Taylor et la compacité du support
de $\widehat{\psi}$, on en déduit (\ref{DL:I:psi}).
\end{dem}
Ainsi, pour une suite $(\lambda_n)_{n \in \N}$ telle que $\ds \lambda_n \limiten \infty$, on définit l'estimateur ponctuel de la densité spectrale
par $\ds\widehat{f}_n(\xi) =
I_{N_n,\psi_{\lambda_n}}\big(\xi^{-1}\big)/
\|\psi\|^2_{L^2} $. On montre que cet estimateur converge vers $f(\xi)$ avec un biais
$\ds \lambda_n^{-1}\times \xi\,f'(\xi) $ et une variance
équivalente à $\ds
\big(\lambda_N
\tau_N^{-1}\big)\times \pi \,  \xi^{-1} f(\xi)
\,\|\widehat{\psi}\|^4_{L^4} /\|\widehat{\psi}\|^4_{L^2}$. Si le
carré du biais est asymptotiquement négligeable devant la
variance, on obtient le TLC suivant ainsi qu'un calcul de l'erreur
quadratique intégrée (MISE):
\begin{Theo}\label{TCL:MISE:dX}
Soient $X$ un processus gaussien défini par
(\ref{repr:harmonizable}) avec une densité spectrale $f$ vérifiant
(F1), (F2) et $\psi$ une ondelette vérifiant  W$(1,0,1/2)$ et $supp
\,\widehat{\psi} \subset [-\Lambda, \Lambda]$. Pour tout réel $\alpha \in ]1/3, 1[$, soit
\ban
\widehat{f}_n(\xi) &:=&  \|\widehat{\psi}\|^{-2}_{L^2}\cdot
\tau_n^{-\alpha} (\xi  \, N_n)^{-1} \sum_{k=1}^{N_n} \left|\int_{\R} e^{i\xi(t-b_k)}\,
\psi\left(\frac {\xi(t-b_k)}{\tau_n^{\alpha}}\right)\, X(t)\,dt\right|^2
\ean
avec la famille de décalages donnée par (\ref{def:bk}).
\begin{description}
\item[i)] Pour toute fréquence $\xi \in ]0,\infty[$, quand $\ds \tau_n\limiten \infty,\quad N_n\limiten \infty\quad
\mathrm{et}\quad \tau_n/N_n
\limiten 0$, on a  \ba \label{TLC2} \tau_n^{(1-\alpha)/2}\big ( \widehat{f}_n(\xi) -  f(\xi) \big ) \limiteloin
\mathcal{N}\Big (0\,,\,4\pi \, \frac {f^2(\xi)} \xi \,  \frac {\|\widehat{\psi}\|^4_{L^4}}
{\|\widehat{\psi}\|^4_{L^2}} \Big ).\ea
\item[ii)] Pour toute bande de fréquences finie $0<\omega_0<\omega_1<\infty$, on a le développement limité du
MISE:
\ban
\E \Big [ \int_{\omega_0}^{\omega_1} \big|\widehat{f}_n(\xi) -  f(\xi)\big|^2\, d\xi \Big ]
    &=&\frac{ 4 \pi}  {\tau_n^{(1-\alpha)}}  \frac{
\|\widehat{\psi}\|^{4}_{L^4}}{ \|\widehat{\psi}\|^{4}_{L^2}} \int_{\omega_0}^{\omega_1}\frac {f^2(\xi)} \xi \,d\xi
\,+\, \tau_n^{-2\alpha}\,\int_{\omega_0}^{\omega_1} \xi^2
f'(\xi)^2\,d\xi\,+\,\mathcal{O}\big(\tau_n^{-2+\alpha}\big).\ean
\end{description}
\end{Theo}
\begin{dem}
La proposition \ref{TCL:ponctuel} et le lemme \ref{lemDL} sont
les deux ingrédients principaux pour obtenir ce théorème.
\end{dem}
\begin{rem} La proposition \ref{TCL:ponctuel}  permet d'estimer des paramètres de la densité spectrale dans un
cadre paramétrique ou semi-paramétrique (par exemple dans le cas du mbf multi-échelle, voir Bardet et Bertrand, 2008). La vitesse de
convergence est alors en $\tau_n^{-1/2}$, où $\tau_n=b_{N_n}-b_1$ correspond à
l'amplitude de la famille de décalages; dans le cadre
non-paramétrique du théorème \ref{TCL:MISE:dX}, la vitesse
de convergence est en $\tau_n^{-1/3-\varepsilon}$ si $\alpha = 1/3 + 2
\varepsilon$ pour un $\varepsilon>0$. Concernant l'estimation du MISE
sur une bande finie $\big[\omega_0, \omega_1\big]$, il n'y a
aucune difficulté à faire tendre la borne supérieure $\omega_1$
vers l'infini; en revanche, on a une explosion à basse fréquence
quand $\omega_0\to 0$.
\end{rem}

\section{Estimation de la densité spectrale pour une trajectoire observée à des instants discrets aléatoires}
\subsection{Autres méthodes existantes}
\subsubsection*{Périodogramme empirique}
L'estimation de la densité spectrale d'un processus gaussien
stationnaire à partir de l'observation d'une trajectoire
échantillonnée à des instants irréguliers
est un vieux
problème en traitement du signal, voir Lii \& Masry (1994) ou Rachdi (2004) et les références qui y sont citées.
Malgré la ressemblance, nous ne travaillons pas sur la même
question. Tout d'abord, nous considérons des processus à
accroissements stationnaires au lieu de processus stationnaires,
mais ce n'est qu'un détail: dans le cas stationnaire, il suffit
de remplacer la formule (\ref{repr:harmonizable}) par $\ds X(t)=
\int_{\R} e^{it\xi} \cdot f^{1/2}(\xi) \, dW(\xi)$.
La vraie différence provient de la modélisation et des
applications sous-jacentes: dans les années 1960, il s'agissait
de choisir une bonne  méthode d'échantillonnage afin d'estimer la
densité spectrale d'un signal analogique; la perspective est
aujourd'hui entièrement différente, on dispose de l'observation
d'un signal à des instants digitaux aléatoires, sans aucune
possibilité de les choisir.

Dans Rachdi (2004) ou Lii \& Masry (1994), on choisit un
échantillonnage poissonien d'intensité moyenne $\beta$ connue,
ainsi que la densité de sa covariance $c(u)$, formules (2.1, 2.2
et 2.3) dans Lii \& Masry (1994), ces quantités sont ensuite
utilisées pour définir la fonction $\gamma(u)$ 
dont
la transformée de Fourier $\Gamma(\lambda)$ intervient dans la
définition de l'estimateur de  la densité spectrale (formule
2.14). Une construction différente, bien que toujours basée sur la
connaissance de la constante $\beta$ et de la fonction $c(u)$ est
proposée dans Rachdi (2004). Cette approche paraît difficilement
transposable à nos données~: en admettant que les temps
d'inter-arrivée soient poissonniens, il faudrait disposer d'une
estimation $\widehat{\beta}$ de $\beta$ et surtout d'une bonne
estimation $\widehat{c}(u)$ de la fonction $c(u)$.

\subsubsection*{Variations quadratiques}
Une autre méthode populaire pour l'estimation du
paramètre de Hurst d'un mbf consiste en l'utilisation des variations quadratiques
généralisées (voir Guyon \& Leon, 1989,  ou Istas \& Lang, 1997).
Pour $(a,b)\in \R_+^2$, soit $Q_X(a,b) := X(b+a) -2 X(b) + X(b-a)$; on estime alors la
variation quadratique d'ordre $2$ (ordre minimal requis pour notre modèle) par $\ds V_N(a) := \sum_{k=1}^N
Q_X(a,b_k)^2$. Le calcul élémentaire
\ban
Q_X(a,b)&=& \int_{\R} e^{ib\xi}\cdot\left(e^{ia\xi}+ e^{-iba\xi}
-2\right)\cdot f^{1/2}(\xi)\, dW(\xi)
\\
&=& 4 \int_{\R} e^{ib\xi}\cdot \sin^2(a\xi/2)\cdot f^{1/2}(\xi)\,
dW(\xi)
\ean
donne une formule de représentation analogue à (\ref{repr:dX}). On
en déduit que le processus $\ds\big(Q_X(a,b)\big)_{b\in
\R}$ est stationnaire, gaussien, centré, de variance
$\ds\mathcal{V}_{X}(a):= 64\, a^{-1} \int_{0}^\infty
\sin^4(2u)\, f(u/a)\, du$, puis que $V_N(a)$ converge vers
$\mathcal{V}_{X}(a)$ selon un TLC analogue à (\ref{TCL1}). Mais l'expression de $\mathcal{V}_{X}(a)$
est telle qu'elle ne permet pas d'estimer la
fonction $f$ à toutes les fréquences (sauf éventuellement en $0$ et en $\infty$ par un choix asymptotique
de $a$).

\subsection{Analyse par ondelette et estimation de la densité spectrale à partir d'une trajectoire observée à des instants discrets aléatoires}
Pour une trajectoire obtenue à des instants discrets $(t_k)_{0\leq k \leq n}$, on définit un
estimateur non-paramétrique de la densité spectrale en remplaçant
les coefficients d'ondelette théoriques par les coefficients
d'ondelette ``discrétisés'' selon la formule (\ref{def:coeff:W:discrets}). On en déduit
\ba
\label{ftilde:non:parametrique}
\widetilde{f}_n(\xi) &=& \xi\, \|\widehat{\psi}\|^{-2}_{L^2}\cdot
\tau_n^{-\alpha} N_n^{-1} \sum_{k=1}^{N_n}\left\{\sum_{i=1}^n \, \left[\int_{t_i}^{t_{i+1}} \cos(\xi(t-b_k))\,
\psi\Big(\frac {\xi(t-b_k}{\tau_n^{\alpha}}\Big)\,dt\right]\times X(t_i)\right\}^2
\ea
Les TLC (\ref{TCL1}) et (\ref{TLC2}) restent vrais, mais les vitesses de
convergence sont plus faibles:
\begin{Theo}\label{tCL:discretisation}
Soit $X$ un processus gaussien défini par
(\ref{repr:harmonizable}) ayant une densité spectrale $f$ vérifiant
(F1) et (F2). On suppose que $\psi$ est une
ondelette vérifiant  W$(1,4,1/2)$ avec $supp
\,\widehat{\psi} \subset [-\Lambda, \Lambda]$, et $X$ est observé selon échantillonnage
vérifiant l'hypothèse (S(s)) avec $3+\big ( 2H-1\, , \, \frac 1 {2H}-\frac 3 2\big )_+ \leq s \leq
\infty$. De plus, s'il existe $\alpha \in ]1/3, 1[$ et
\begin{itemize}
\item si $s=\infty$, $n \, \delta_n^{2+H}\limiten 0$ et $\tau_n \, n^{-(1+H)/(2+H)} \limiten 0$,
\item si $s<\infty$, $n\,   \delta_n^{2+H-\frac {(H+1)^2}{H+s}}\hspace{-3mm}\limiten 0$ et $\tau_n \, n^{-(s(H+1)+2\alpha)/(s(H+2)-1+\alpha(1-H))} \limiten 0$,
\end{itemize}
\begin{description}
\item[i)] Le TLC (\ref{TLC2}) est vérifié également par $\widetilde{f}_n(\xi)$;
\item[ii)] Pour toute bande de fréquences finie $0<\omega_0<\omega_1<\infty$, on a le développement limité du
MISE:
\ban
\E \Big [ \int_{\omega_0}^{\omega_1} \big|\widehat{f}_n(\xi) -  f(\xi)\big|^2\, d\xi \Big ]
    &=&\frac{ 4 \pi}  {\tau_n^{(1-\alpha)}}  \frac{
\|\widehat{\psi}\|^{4}_{L^4}}{ \|\widehat{\psi}\|^{4}_{L^2}} \int_{\omega_0}^{\omega_1}\frac {f^2(\xi)} \xi \,d\xi
\,+\, \tau_n^{-2\alpha}\,\int_{\omega_0}^{\omega_1} \xi^2
f'(\xi)^2\,d\xi\,+\,\mathcal{O}\big(\tau_n^{-2+\alpha)}\big).\ean\end{description}
\end{Theo}
\section{Applications numériques}
\subsection*{Résultats de simulations pour des mouvements browniens fractionnaires}
La densité spectrale d'un mouvement brownien fractionnaire de paramètres $(H,\sigma^2)$ s'exprime sous la forme $f_{H,\sigma^2}(\xi)=\pi^{-1} H \, \Gamma(2H)\, \sin(\pi H) \, \sigma^2\, |\xi|^{-1-2H}$. Nous avons simulé
des trajectoires de différentes longueurs de mbf pour différentes valeurs de $H$, avec
$\sigma^2= \pi \big (H \, \Gamma(2H)\, \sin(\pi H)\big )^{-1}$, de telle manière que la densité spectrale soit
$ |\xi|^{-1-2H}$.
Il est à noter que dans une telle simulation le choix de $\alpha$ s'avère crucial:
si $\alpha$ est choisi trop proche de la valeur frontière $1/3$, l'algorithme d'estimation n'est pas stable;
de même si $\alpha$ est trop proche de $1$, l'estimateur ne converge pas assez vite et sa variance
est très importante. De multiples simulations nous ont montré que plus $N$ est grand, plus $\alpha$ pouvait être
choisi proche de $1/3$. La Figure~1 donne un exemple d'une telle estimation pour $N=20000$, $H=0.2$ et $\alpha=0.6$.
\par
\vspace{1.5cm}
\begin{figure}[h]
\[
\epsfxsize 8cm \epsfysize 6cm \epsfbox{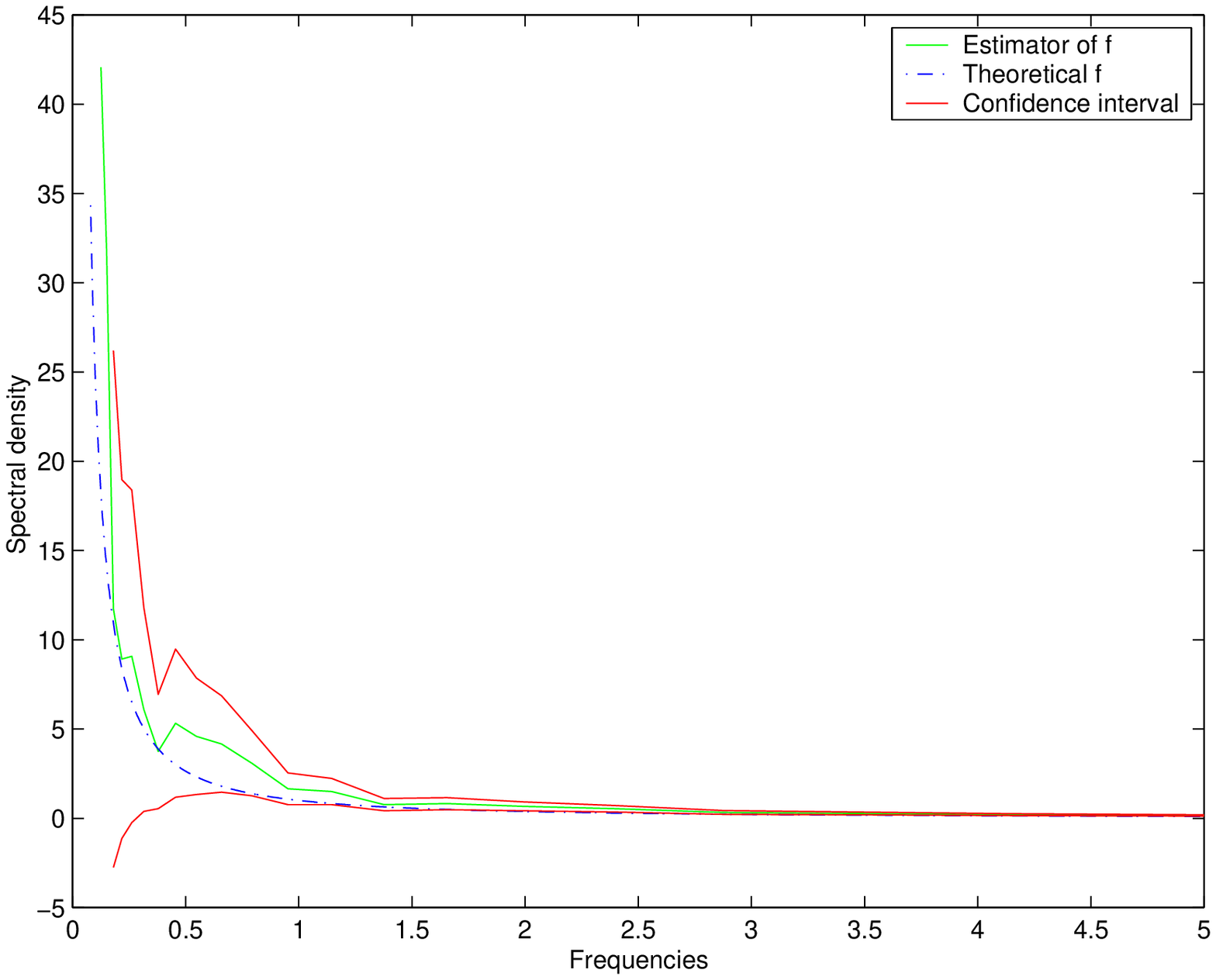}
\hspace*{1.5 cm} \epsfxsize 8cm \epsfysize 6cm
\epsfbox{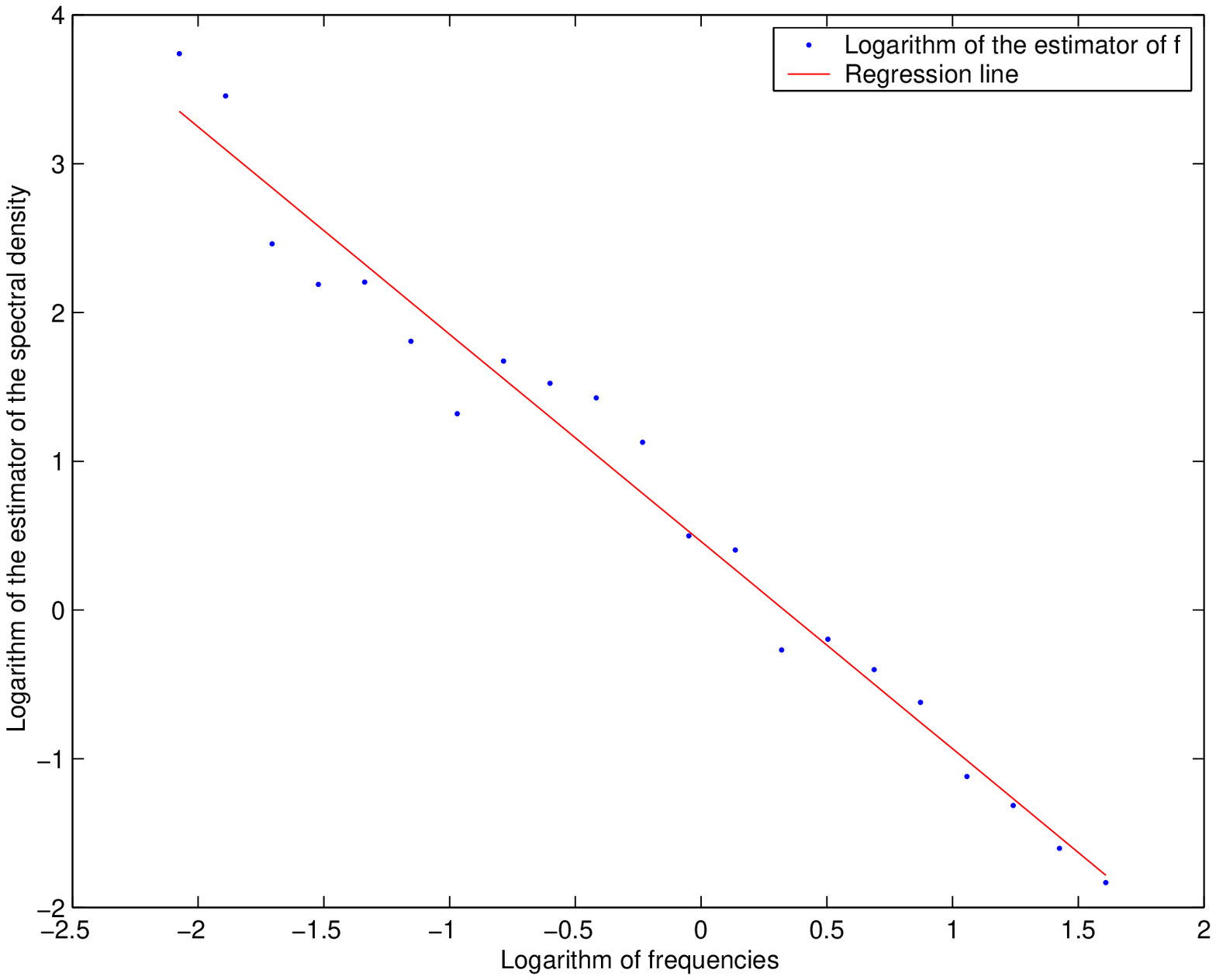}
\]
\label{figFBM}
\caption{\it Estimation de la densité spectrale (à gauche) et de son logarithme (à droite) pour une trajectoire de  mbf telle
que $H=0.2$ et $N=20000$ échantillonné suivant un pas suivant une loi exponentielle de paramètre $1$.}
\end{figure}
\subsection*{Application aux fréquences cardiaques instantanées pendant un marathon}
Les données que nous considérons ici ont été obtenues par le Laboratoire INSERM LEPHE, dirigé
par la professeure V. Billat (Université d'Evry, France). Il s'agit des mesures (parmi d'autres mesures obtenues) des
durées entre deux battements de coeur de $50$ marathoniens de bon niveau courant le Marathon de Paris 2004.
Ces données ont également été étudiées dans Kammoun {\it et al.} (2007), en considérant la suite
des fréquences cardiaques instantanées comme une série chronologique classique (à pas régulier).
Ici, nous considérons l'évolution de la fréquence cardiaque instantanée en fonction du temps réellement
écoulé (suivant un pas de temps irrégulier correspondant aux durées entre deux battements successifs).
La Figure~2 présente un exemple sur un athlète (ayant donc couru ce marathon en 2h45mn12s)~:
\par
\vspace{1cm}
\begin{figure}[h]
\[
\epsfxsize 8cm \epsfysize 6cm \epsfbox{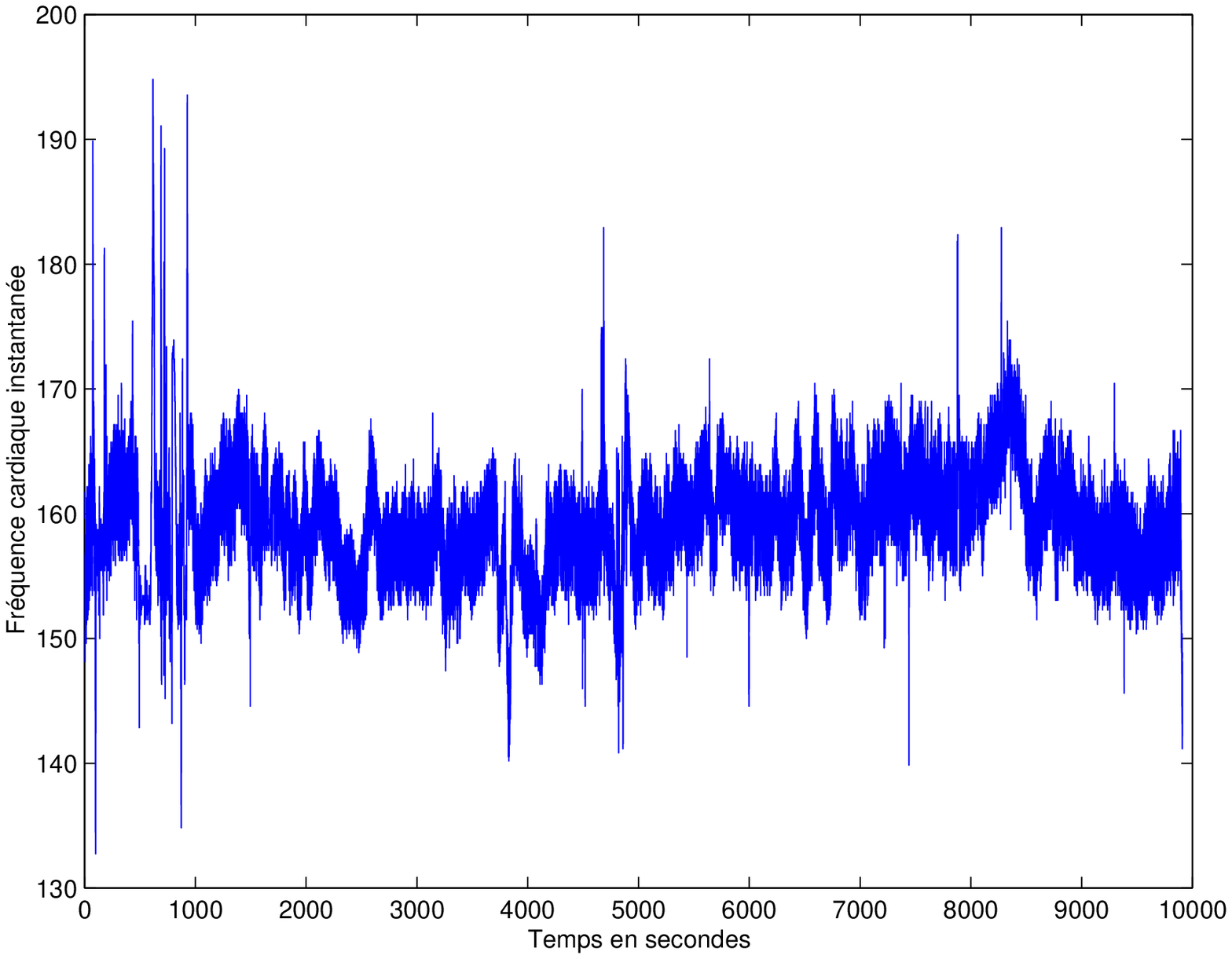}
\hspace*{1.5 cm} \epsfxsize 8cm \epsfysize 6cm
\epsfbox{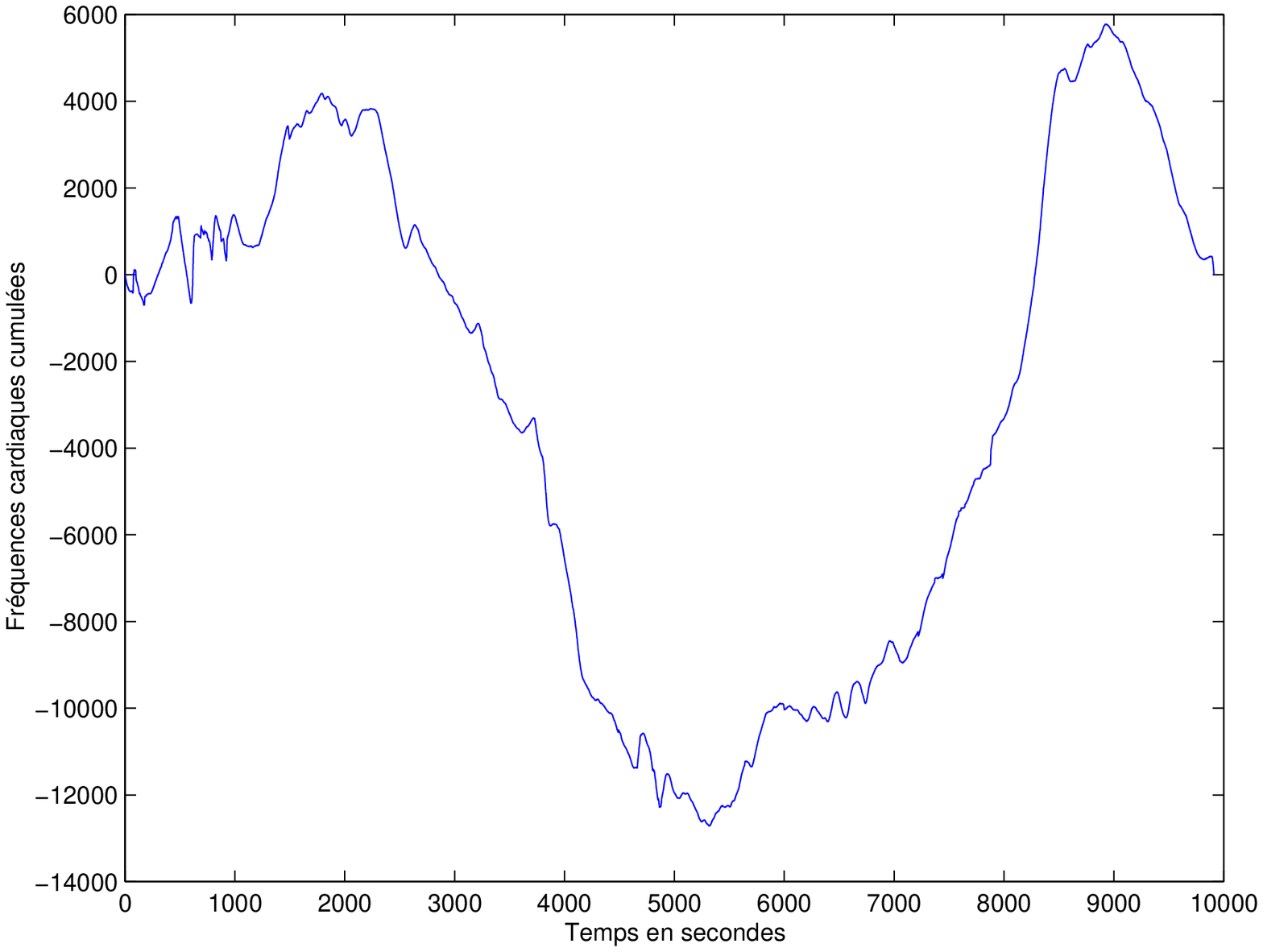}
\]
\label{Athlete1}
\caption{ \it Fréquences cardiaques instantanées (à gauche) et cumulées (à droite) d'un marathonien en fonction du temps de course}
\end{figure}
\noindent Nous avons appliqué l'estimateur de la densité spectrale $\widetilde{f}_n(\xi) $ aux
des fréquences cardiaques instantanées de ce coureur (fréquences cumulées qui semblent satisfaire les hypothèses du modèle de processus gaussien à accroissements stationnaires considérées ici). Cependant, comme
cela avait déjà été proposé dans Kammoun {\it et al.} (2007), on a préféré découper
cette série en trois sous-séries (obtenues à partir d'un algorithme de détection de ruptures en moyenne
et variance) représentant le début ($4802$ données), le milieu ($19708$ données) et la fin de course ($2590$ données)  afin de tester une éventuelle évolution de l'indice d'auto-similarité du signal fréquence cardiaque en fonction de la fatigue.
Notons que le pas d'échantillonnage ici suit approximativement une loi gaussienne (donc $S(s)$ est vérifié pour tout $s>0$) de moyenne $0.376s$ et d'écart-type $0.01s$. Pour appliquer les résultats précédents, on considère la base de temps en minutes et ainsi $\delta_n \simeq 0.006 mn$. La Figure~3 ci-dessous représente le logarithme de la densité spectrale en fonction du logarithme des fréquences pour la partie de fin de course.
\par
\vspace{2cm}
\begin{figure}[h]
\[
\epsfxsize 10cm \epsfysize 8cm \epsfbox{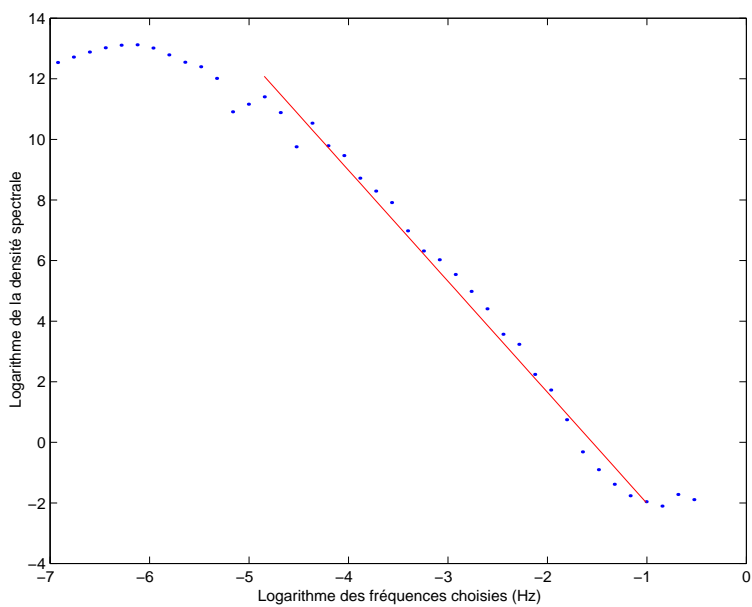}
\]
\label{Athlete_fin}
\caption{\it Logarithme de la densité spectrale en fonction du logarithme des fréquences (fin de course)}
\end{figure}
On retrouve dans la Figure~3 la même propriété évoquée dans Kammoun {\it et al.} (2007): il
existe une phase de linéarité dans ce graphe. Ceci avait permis la modélisation de ces données par
un mouvement brownien localement fractionnaire, qui est une généralisation du MBF (même densité
spectrale dans une bande de fréquence) et le coefficient directeur de la droite fournit une estimation
du paramètre de fractalité locale $H$. On trouve, de la même manière:
$$
\widehat H_{D}\simeq 1.12~\mbox{(début de course)},\quad \widehat H_{M}\simeq 1.15~\mbox{(milieu de course)} ~~\mbox{et}~~ \widehat H_{D}\simeq 1.38~\mbox{(fin de course)}.
$$
On constate une augmentation de la valeur de ce paramètre en fin de course, ce qui donne une
intéressante caractérisation de l'apparition d'une fatigue du coeur (voir Kammoun {\it et al.}, 2007). Ce résultat devra être étayé par des études ultérieures afin par exemple d'expliquer  l'augmentation  depuis une décennie de l'occurence de mort subite au cours des marathons.
\par
\bigskip
\noindent{\bf Remerciements~:} Nous tenons à remercier Mario Wschebor (Université de Montevideo) et
Antoine Ayache (Université de Lille I)  dont les questions et conseils ont contribué à la mise au point de l'algorithme et Imen Kammoun( Université de Paris I) pour le pré-traitement des données. Cependant, toutes les erreurs éventuelles sont de notre seule responsabilité.



\begin{thebibliography}{99}

\bibitem{Abry} Abry, P., Flandrin, P., Taqqu, M.S. and Veitch, D. (2002).
Self-similarity and long-range dependence through the wavelet
lens, in  {\it Long-range Dependence: Theory and Applications}, P.
Doukhan, G. Oppenheim and M.S. Taqqu editors, Birkh{\"a}user.

\bibitem{YAS:PM:08}  Aït-Sahalia and Y. \& Mykland, P.A. (2008). An analysis of Hansen-Scheinkman moment estimators for
discretely and randomly sampled diffusions,  {\em Journal of
Econometrics}, In Press.


\item Bardet, J.M. and Bertrand, P. (2007), "Identification of the
multiscale fractional Brownian motion with biomechanical
applications", {\it Journal of Time Series Analysis}, 28, p. 1-52.


\bibitem{BB:08} Bardet, J.M. and Bertrand, P.R. (2008), " Wavelet Analysis along Random Sampling for Gaussian
Processes". Preprint.


\bibitem{Begyn05} Begyn, A. (2005). Quadratic Variations along
Irregular Subdivisions for Gaussian Processes. {\em Electronic
Journal of Probability}, 10, p-691-717.

\bibitem{Ben}
Bertrand, P. and Bardet, J.M. (2001), " Some generalization of fractional
Brownian motion
 and Control", in {\it Optimal Control and Partial Differential
 Equations}, J.L.~Menaldi, E.~Rofman and A.~Sulem editors, p.221-230,  IOS
 Press .

\bibitem{Blanke:SPA} D. Blanke, and Vial, C. (2008), "Assessing the
number of mean-square derivative of a Gaussian process", to
appears in {\em Stoch. Proc. Applications}.





\bibitem{GVH:02} Gao, J.; Anh, V.; Heyde, C. (2002)
Statistical estimation of non-stationary Gaussian processes with long-range dependence and intermittency.
{\it Stochastic Process. Appl.} {\bf  99}, no. 2, 295--321

\bibitem{Guyon:Leon:89}
Guyon, X., and Leon, J.R.(1989). Convergence en loi des
H-variations d'un processus gaussien fractionnaire, {\it Ann.
Inst. H. Poincar\'e}, 25,  265-282.

\bibitem{YH:05} Hayashi, T. and Yoshida, N. (2005).  On covariance estimation of
non-synchronously observed diffusion processes. {\em Bernoulli},
Volume 11, Number 2, p.359-379.

\bibitem{is-lan} Istas, J. and Lang, G. (1997). Quadratic variations
and estimation of the local H\"older index of a Gaussian process.
{\em Ann. Inst. Poincar\'e}, 33, 407-436.


\bibitem{KBB:07} Kammoun, I., Billat, V and Bardet, J.M. (2007). Comparison of DFA vs wavelet analysis for estimation of regularity of HR series during the marathon. Preprint available on {\tt http://hal.archives-ouvertes.fr/}.


\bibitem{Lii:M:94}
 Lii, K.S. and Masry, E. (1994) Spectral estimation of continuous-time stationary processes from random sampling.
 {\it Stochastic Process. Appl.} {\bf 52}, no. 1, 39--64.


\bibitem{Radch:04}
 Rachdi, M.(2004).
Strong consistency with rates of spectral estimation of
continuous-time processes: from periodic and Poisson sampling
schemes. {\it J. Nonparametr. Stat. } {\bf 16} , no. 3-4,
349--364.

\bibitem{TaSa:1994} Samorodnitsky, G. and Taqqu M.S. (1994), {\it Stable
non-Gaussian Random Processes}, Chapman and Hall.



\bibitem{Ya:55} Yaglom, A.M. (1958). Correlation theory of processes with stationary random increments of order $n$. {\it Trans. A.M.S.} 8, p.87-141.


\end{thebibliography}
\end{document}